\documentclass[12pt]{amsart}

\newtheorem{THM}{Theorem}[section]
\newtheorem{LMA}[THM]{Lemma}
\newtheorem{PROP}[THM]{Proposition}
\newtheorem{CORO}[THM]{Corollary}
\newtheorem{CONJ}[THM]{Conjecture}
\newtheorem{EG}[THM]{Example}
\newtheorem{PROB}[THM]{Problem}

\numberwithin{equation}{section}

\usepackage{amsmath}
\usepackage{amssymb}
\usepackage{euscript}

\newcommand{\showon}{\begin{eqnarray*}}
\newcommand{\showoff}{\end{eqnarray*}}

\newcommand{\NN}{\mathbf{N}}
\newcommand{\RR}{\mathbf{R}}
\newcommand{\ZZ}{\mathbf{Z}}
\newcommand{\CC}{\mathbf{C}}

\newcommand{\C}{\EuScript{C}}
\newcommand{\D}{\EuScript{D}}
\newcommand{\K}{\EuScript{K}}
\newcommand{\T}{\EuScript{T}}
\newcommand{\G}{\EuScript{G}}
\renewcommand{\S}{\EuScript{S}}
\renewcommand{\a}{\mathbf{a}}
\renewcommand{\b}{\mathbf{b}}
\renewcommand{\c}{\mathbf{c}}
\newcommand{\h}{\mathbf{h}}
\newcommand{\m}{\mathbf{m}}
\newcommand{\n}{\mathbf{n}}
\newcommand{\q}{\mathbf{q}}
\renewcommand{\u}{\mathbf{u}}
\renewcommand{\v}{\mathbf{v}}
\newcommand{\w}{\mathbf{w}}
\newcommand{\x}{\mathbf{x}}
\newcommand{\y}{\mathbf{y}}
\newcommand{\z}{\mathbf{z}}

\newcommand{\one}{\boldsymbol{1}}
\newcommand{\zero}{\boldsymbol{0}}
\newcommand{\drop}{\smallsetminus}
\newcommand{\goesto}{\rightarrow}

\newcommand{\diag}{\mathrm{diag}}
\newcommand{\aut}{\mathrm{Aut}}
\newcommand{\snf}{\mathrm{snf}}

\newcommand{\bank}{\$}

\begin{document}

\title[The critical group of a directed graph]{The critical group of a
directed graph}
\author{David G. Wagner}
\address{Department of Combinatorics and Optimization\\
University of Waterloo\\
Waterloo, Ontario, Canada\ \ N2L 3G1}
\email{\texttt{dgwagner@math.uwaterloo.ca}}
\thanks{Research supported by the Natural
Sciences and Engineering Research Council of Canada under
operating grant OGP0105392.}

\begin{abstract}
For a finite, directed graph $G=(V,E)$ we define the critical group $\K(G)$
to be the cokernel of the transpose of the Laplacian matrix of $G$ acting on 
$\ZZ^{V}$, and $K(G)$ to be its torsion subgroup.  This generalizes the
case of undirected graphs studied by Bacher, de la Harpe and Nagnibeda,
and by Biggs.  We prove a variety of results about these critical groups,
among which are:  that $\K(G/\pi)$ is a subgroup of $\K(G)$ when $\pi$ is
an equitable partition and $G$ is strongly connected;  that $K(G)$ depends
only on the graphic matroid of $G$ when $G$ is undirected;  that
there is no `natural' bijection between spanning trees of $G$ and $K(G)$ when
$G$ is undirected, even though these sets are equicardinal; and that
the `dollar game' of Biggs can be generalized slightly to provide a
combinatorial interpretation for the elements of $K(G)$ when $G$ is
strongly connected.
\end{abstract}
\maketitle

\section{Introduction.}

We use the word \emph{graph} to refer to a finite, possibly
directed, multigraph.  If a graph is undirected then we consider each 
of its edges to represent a pair of directed edges with opposite
orientations.  (In particular, an undirected loop represents two
directed loops.)  The \emph{adjacency matrix} $A(G)$ of a graph 
$G=(V,E)$ is indexed by $V\times V$, with $vw$--entry $A_{vw}$ being
the number of directed edges of $G$ with initial vertex $v$ and terminal
vertex $w$.  The matrix $\Delta(G)$ is indexed by $V\times V$, with
diagonal entry $\Delta_{vv}:=\sum_{w\in V}A_{vw}$
being the outdegree of the vertex $v$, and with zero off-diagonal entries.
The \emph{Laplacian matrix} of $G$ is $Q(G):=\Delta(G)-A(G)$.

The (\emph{full}) \emph{critical group} $\K(G)$ of $G=(V,E)$
is the cokernel of the transpose of its Laplacian matrix acting
on $\ZZ^{V}$;  that is,
$$\K(G):=\ZZ^{V}/Q^{\dagger}(G)\ZZ^{V}.$$
This is a finitely generated abelian group.
Also, we define the \emph{reduced critical group} $K(G)$ of $G$ to
be the torsion subgroup of $\K(G)$;  that is, the subgroup of $\K(G)$
consisting of all elements of finite order.  For a connected, undirected
graph $G$, this $K(G)$ is the `Jacobian group' defined by Bacher,
de la Harpe, and Nagnibeda \cite{BHN} and studied further by Biggs
\cite{Bi2,Bi3,Bi4} under the term `critical group'.  Our goal here is to
investigate the relationship between the combinatorial structure of
graphs and the algebraic structure of their critical
groups as generally as possible.  Since the Laplacian
matrix is insensitive to loops in a graph, we might as well
restrict attention to loopless graphs;  it is only in Section 9, 
however, that we really require this restriction.

Section 2 contains some preliminary observations.  In Section 3,
we determine the rank of $\K(G)$ combinatorially.  In Section 4,
we consider the minimal number of generators of $K(G)$;
this is a much more difficult invariant of $G$ than the rank of
$\K(G)$, and we obtain only a weak upper bound for it.
In Sections 5 and 6, we show that $\K(G/\pi)$ is (isomorphic
to) a subgroup of
$\K(G)$ when $\pi$ is an equitable partition of a strongly
connected graph $G$.  In Section 7, we prove some isomorphism
theorems, with the consequence that for an undirected graph $G$,
the reduced critical group $K(G)$  depends only on the
graphic matroid of $G$.  In Section 8, we show that, for a 
connected undirected graph $G$, there is no `natural' bijection
between $K(G)$ and the set of spanning trees of $G$, even though
these sets are equicardinal.  In Section 9, we revisit the `dollar
game' of Biggs \cite{Bi2,Bi3,Bi4} in the more general setting of
strongly connected graphs.  The theory is almost the same as for
undirected graphs, with  one interesting new complication.
Throughout the paper, we indicate various conjectures and open
problems.

\section{Preliminaries.}

Let $G$ be a graph with $n(G)$ vertices.  By the structure theorem for
finitely generated abelian groups, there are nonnegative integers
$g_{1}$,\ldots, $g_{n}$ such that $g_{i}$ divides $g_{i+1}$ for each
$1\leq i\leq n-1$ and
$$\K(G)\simeq (\ZZ/g_{1}\ZZ)\oplus(\ZZ/g_{2}\ZZ)\oplus\cdots\oplus
(\ZZ/g_{n}\ZZ).$$
(Of course, $\ZZ/1\ZZ=0$ is the trivial group and $\ZZ/0\ZZ=\ZZ$.)
These integers are computed by reducing $Q(G)$ to its Smith
normal form.   Say that two $n$-by-$n$ integer matrices $M$ and $M'$
are \emph{equivalent}, denoted by $M\approx M'$, if and only if there
exist $n$-by-$n$ integer matrices $L$ and $N$ with determinant $\pm 1$
such that $LMN=M'$.  That is, $M$ can be transformed to $M'$ by 
applying elementary row and column operations which are invertible over
the integers.  It is not difficult to verify that every $n$-by-$n$
integer matrix $M$ is equivalent to a diagonal matrix $\diag(g_{1},g_{2},
\ldots,g_{n})$ of nonnegative integers such that $g_{i}$ divides
$g_{i+1}$ for each $1\leq i\leq n-1$, and that this matrix $\snf(M)$ is
determined uniquely by $M$.  This is the \emph{Smith normal form} of 
$M$.

We define the \emph{dual critical group} of $G$ to be
$$\K^{*}(G):=\ZZ^{V}/Q(G)\ZZ^{V},$$
the cokernel of the Laplacian of $G$ acting on $\ZZ^{V}$.
It is clear that for any square integer matrix $M$,
$\snf(M^{\dagger})=\snf(M)$, from which Proposition 2.1 follows.
\begin{PROP}
For any graph $G$, $\K^{*}(G)$ is isomorphic to
$\K(G).$
\end{PROP}
Of course, $G$ is undirected if and only if $Q^{\dagger}=Q$,
in which case $\K^{*}(G)=\K(G)$.  In general, however,
the isomorphism in Proposition 2.1 depends on a
choice of matrices $L$ and $N$ such that $LQN=\snf(Q)=N^{\dagger}
Q^{\dagger}L^{\dagger}$, and hence is not natural.  These
dual critical groups will be useful in Section 6.

If $G$ is a graph with weak components $G_{1},\ldots,G_{c}$, then
$Q(G)=Q(G_{1})\oplus\cdots\oplus Q(G_{c})$ is block-diagonal, from
which Proposition $2.2$ follows.
\begin{PROP}
If $G$ is a graph with weak components $G_{1},\ldots,G_{c}$, then
$$\K(G)=\K(G_{1})\oplus\cdots\oplus\K(G_{c}).$$
\end{PROP}

Proposition $2.3$ is essentially the Matrix-Tree Theorem; see Theorem 
6.3 of Biggs \cite{Bi1} or Theorem 7.3 of Biggs \cite{Bi3}.
\begin{PROP}
If $G$ is connected and undirected then the order of $K(G)$
is $\kappa(G)$, the number of spanning trees of $G$.
\end{PROP}

We use the slightly odd but convenient notations $M_{vw}$
for the  $(v,w)$-entry of a $V$-by-$V$ matrix $M$, but $x(v)$
for the $v$-th entry of a $V$-indexed vector $\x$.  Also,
$\one$ denotes the all-ones vector, and $\zero$ denotes the
zero vector.

\section{The torsion-free part of $\K(G)$.}

For each natural number $g$, let $\mu_{g}(G)$ denote the multiplicity
with which $g$ occurs on the diagonal of $\snf(Q^{\dagger}(G))$.  Thus,
$\mu_{0}(G)$ is the rank of $\K(G)$, so that
$\K(G)\simeq K(G)\oplus\ZZ^{\mu_{0}(G)}.$
In this section we determine the combinatorial meaning of $\mu_{0}(G)$
for any graph $G$. 

Since $$(\ZZ/g\ZZ)\otimes\RR=\left\{\begin{array}{ll}
\RR &\ \mathrm{if}\ g=0,\\
0 &\ \mathrm{if}\ g\neq 0,\end{array}\right.$$
for all nonnegative integers $g$, and since tensor product
distributes across direct sums, we see that
$$\mu_{0}(G) = \dim_{\RR}\K(G)\otimes\RR=\dim_{\RR}\RR^{V}/Q^{\dagger}\RR^{V}
=\dim_{\RR}\ker(Q^{\dagger}).$$
The results of this section, determining $\dim_{\RR}\ker(Q^{\dagger})$
combinatorially, are well-known, but we repeat the short proofs for
completeness and the readers' convenience.

\begin{LMA}  Let $G$ be a strongly connected graph.\\
\textup{(a)} The kernel of $Q$ acting on $\RR^{V}$ is $\RR\one$, the 
span of the all-ones vector.\\
\textup{(b)} If $T$ is a diagonal matrix with nonnegative real
entries, then either $T=O$ or $Q+T$ is invertible over $\RR$.
\end{LMA}
\begin{proof}
We prove (a) and (b) together by showing that for $T$ as in part (b)
and $\z\neq\zero$, if $(Q+T)\z=\zero$ then $T=O$ and $\z=c\one$
for some $c\in\RR$.  With these hypotheses, choose any vertex $v\in V$
such that $|z(v)|>0$ is maximum.  Then, since $(\Delta+T)\z=A\z$ we see
that
$$(\Delta_{vv}+T_{vv})z(v)=\sum_{w\in V}A_{vw}z(w).$$
Since there are $\Delta_{vv}$ terms on the right side (considering
$A_{vw}$ as the multiplicity of the term $z(w)$) and each of these
has absolute value at most $|z(v)|$, it follows that $T_{vv}=0$
and $z(w)=z(v)$ for all $w\in V$ such that $A_{vw}\neq 0$.  Now, we
may repeat this argument with any such vertex $w$ in place of $v$,
\emph{et cetera}.  Since $G$ is strongly connected, it follows that
$T=O$ and $\z=z(v)\one$, completing the proof.
\end{proof}

\begin{PROP}  Let $G$ be a strongly connected graph.  There is a 
unique vector $\h\in\RR^{V}$ such that $Q^{\dagger}\h=\zero$, the
entries of $\h$ are positive integers, and $\mathrm{gcd}\{
h(v):\ v\in V\}=1$.  Moreover, $\ker(Q^{\dagger})=\RR\h$.
\end{PROP}
\begin{proof}
If $G$ has a single vertex then the result is trivial, so assume
$n(G)\geq 2$.  Now, since $G$ is strongly connected, the matrix
$\Delta$ is invertible.
Since $\ker(Q)$ is one-dimensional, $\ker(Q^{\dagger})$ is also
one-dimensional;  hence there is a unique (nonzero) vector $\z$
such that $Q^{\dagger}\z=\zero$ and $\one^{\dagger}\Delta\z=1$.
Then 
$(\Delta\z)^{\dagger}\Delta^{-1}A=\z^{\dagger}A=(\Delta\z)^{\dagger}$,
so $\Delta\z$ is the stationary distribution of the Markov
chain represented by the stochastic matrix $\Delta^{-1}A$.
Since $G$ is strongly connected, every state of this Markov chain
is recurrent, so every entry of $\z$ is positive.  Since $\z$ solves
the system $Q^{\dagger}\z=\zero$, which has integer coefficients,
every entry of $\z$ is rational.  Finally, there is a unique
positive integer multiple of $\z$ which gives the vector $\h$
with the desired properties.
\end{proof}
For a strongly connected graph $G$, the vector $\h$ defined
in Proposition 3.2 will be significant for several results in what follows.
We refer to $h(v)$ as the \emph{activity} of the vertex $v\in V$,
for reasons which will be seen in Section 9.

For $S$ a strong component of $G$, if $\x\in\RR^{V}$ then let $\x|_{S}$
be the restriction of $\x$ to $V(S)$, and if $M$ is a $V$-by-$V$ 
matrix then let $M|_{S}$ denote the submatrix of $M$ indexed by
rows and columns from $V(S)$.

\begin{LMA} For any graph $G$, if
$\z\in\ker(Q^{\dagger})$ and $S$ is a non-terminal strong component
of $G$, then $\z|_{S}=\zero$.
\end{LMA}
\begin{proof}
Let $S_{1},\ldots,S_{c}$ be a list of the strong components of $G$
such that if there is a directed edge from $v\in S_{i}$ to $w\in 
S_{j}$, then $i\leq j$.  We prove the claim by induction on $1\leq j\leq c$.
For the basis of induction ($j=1$), and for the induction step
($j\geq 2$), we may assume that $\z|_{S_{i}}=\zero$ for all $1\leq 
i<j$ such that $S_{i}$ is a non-terminal strong component of $G$.
If $S_{j}$ is a terminal strong component of $G$ then there is
nothing to prove. Otherwise, there is at least
one directed edge with initial vertex in $S_{j}$ and terminal vertex
not in $S_{j}$.  Therefore, $Q(G)|_{S_{j}}=Q(S_{j})+T$ for some
nonzero diagonal matrix $T$ of nonnegative integers.  Now, since
$\z|_{S_{i}}=\zero$ for all $1\leq i<j$ such that $S_{i}$  
non-terminal,
$$\zero=(Q^{\dagger}(G)\z)|_{S_{j}}=(Q^{\dagger}(S_{j})+T)(\z|_{S_{j}}),$$
and by Lemma 3.1(b) we conclude that $\z|_{S_{j}}=\zero$.
\end{proof}

\begin{THM}  Let $G$ be any graph,
and let the terminal strong components of $G$ be $S_{1},\ldots,S_{t}$.
Then
$$\ker(Q^{\dagger}(G))\simeq\ker(Q^{\dagger}(S_{1}))\oplus
\cdots\oplus\ker(Q^{\dagger}(S_{t})).$$
\end{THM}
\begin{proof}
Let $\z\in\ker(Q^{\dagger}(G))$.  By Lemma 3.3, $\z|_{S}=\zero$ if
$S$ is not a terminal strong component of $G$.  Hence, if $S$ is a 
terminal strong component of $G$, then
$$\zero=(Q^{\dagger}\z)|_{S}=Q^{\dagger}(S)(\z|_{S}),$$
so that $\z|_{S}\in\ker(Q^{\dagger}(S)).$  Conversely, if
$\z_{i}\in\ker(Q^{\dagger}(S_{i}))$ for each terminal
strong component of $G$ then the $\z\in\RR^{V}$ defined by
$$\z|_{S}:=\left\{\begin{array}{ll}
\z_{i} & \mathrm{if}\ S=S_{i}\ \mathrm{for\ some}\ 1\leq i\leq t ,\\
\zero & \mathrm{otherwise},\end{array}\right.$$
is in the kernel of $Q^{\dagger}(G)$.
\end{proof}

\begin{CORO}  For any graph $G$, $\mu_{0}(G)$ is the number of
terminal strong components of $G$.
\end{CORO}

\section{The minimal number of generators of $K(G)$.}

Let $\nu(G)$ denote the minimal number of generators of the reduced
critical group of $G$, so $K(G)$ is the direct sum of $\nu(G)$ nontrivial
finite cyclic groups.  This is a rather difficult numerical invariant
of $G$, as it depends on the arithmetic properties of $Q(G)$.
Clearly, $\mu_{1}(G)+\nu(G)+\mu_{0}(G)=n(G)$.  We give an upper bound
on $\nu(G)$ by proving a combinatorial lower bound on $\mu_{1}(G)$.  
It must be admitted, however, that this bound is generally quite
weak.

Let $G$ be a weakly connected graph, with Laplacian
matrix $Q$.  A \emph{reduction sequence} in $G$ is a sequence
$(v_{1},w_{1}),\ldots,(v_{s},w_{s})$ of pairs of vertices such that:\\
$\bullet$\ the vertices $v_{1},\ldots,v_{s}$ are pairwise distinct,\\
$\bullet$\ the vertices $w_{1},\ldots,w_{s}$ are pairwise distinct,\\
$\bullet$\ for each $1\leq i\leq s$,\ $Q_{v_{i}w_{i}}\in\{-1,1\}$, and\\
$\bullet$\ for each $2\leq j\leq s$,\ either $Q_{v_{i}w_{j}}=0$ for
all $1\leq i<j$, or $Q_{v_{j}w_{i}}=0$ for all $1\leq i<j$.\\
Let $\sigma(G)$ denote the maximum length of a reduction sequence in $G$.

\begin{PROP}
If $G$ is a weakly connected graph, then
$\mu_{1}(G)\geq\sigma(G)$, and hence $\nu(G)\leq n(G)-\mu_{0}(G)-\sigma(G)$.
\end{PROP}
\begin{proof}
Let $(v_{1},w_{1}),\ldots,(v_{s},w_{s})$ be a reduction sequence in
$G$, and let $Q$ be the Laplacian matrix of $G$.

The first claim is that we may apply elementary row and column
operations to $Q$, involving only rows $v_{1},\ldots,v_{s}$ and columns
$w_{1},\ldots,w_{s}$, so that these rows and columns of the resulting
matrix induce an $s$-by-$s$ identity matrix.  We prove this by induction
on $s$, the basis $s=1$ being clear.  For the induction step, since
$(v_{1},w_{1}),\ldots,(v_{s-1},w_{s-1})$ is a reduction sequence of
length $s-1$, the induction hypothesis gives elementary
operations involving only rows $v_{1},\ldots,v_{s-1}$ and columns
$w_{1},\ldots,w_{s-1}$, which, when applied to $Q$, result in a matrix
$Q'$ in which these rows and columns induce an $(s-1)$-square identity
matrix.  By the last condition defining a reduction sequence, we have
either $Q_{v_{i}w_{s}}=0$ for all $1\leq i<s$, or $Q_{v_{s}w_{i}}=0$
for all $1\leq i<s$.  Examining the way in which $Q'$ was obtained from
$Q$, we see that either $Q'_{v_{i}w_{s}}=0$ for all $1\leq i<s$, or
$Q'_{v_{s}w_{i}}=0$ for all $1\leq i<s$.  Now, elementary row or
column operations can be used to cancel any nonzero entries
$Q'_{v_{i}w_{s}}\neq 0$ or $Q'_{v_{s}w_{i}}\neq 0$ with $1\leq i<s$,
and the value of $Q'_{v_{s}w_{s}}$ is left unchanged.  Finally, multiplying row
$v_{s}$ by $-1$ if necessary produces the $s$-by-$s$ identity submatrix,
as claimed.

Having produced this $s$-by-$s$ identity submatrix, we use
elementary column operations to zero out all
entries in rows $v_{1},\ldots,v_{s}$ except for the $1$s in the
$(v_{i},w_{i})$ positions $(1\leq i\leq s)$.  Then we use
elementary row operations to zero out all
entries in columns $w_{1},\ldots,w_{s}$ except for the $1$s in the
$(v_{i},w_{i})$ positions $(1\leq i\leq s)$.  The result is a matrix,
equivalent to $Q$, which is also equivalent to a matrix with the
block structure $I_{s}\oplus M$ for some $(n-s)$-square matrix $M$.
Therefore, $1$ occurs at least $s$ times in the Smith normal form
of $Q$.  Considering a reduction sequence of maximum length $s=\sigma(G)$
completes the proof.
\end{proof}

The bound of Proposition 4.1 is likely to be very far from the
true value of $\nu(G)$, since it makes no use of the arithmetic
structure of $Q(G)$.  With this in mind, here are
a few conjectures.  Let $\G(n,p)$ denote a random
simple, undirected graph with $n$ vertices and edge-probability
$0\leq p\leq 1$.  As is is well-known (see Theorem 4.3.1 of
Palmer \cite{Pa}),
if $p(n)>c\log(n)/n$ with $c>1$ then, as $n\goesto\infty$, the
probability that $\G(n,p)$ is connected converges to $1$.

\begin{CONJ}
If $c>1$ and $c\log(n)/n<p(n)<1-o(\log(n)/n)$ then, as
$n\goesto\infty$, the probability that $K(\G(n,p))$ is
cyclic converges to $1$.
\end{CONJ}
That is, the conjecture is that almost every connected undirected
simple graph has a cyclic reduced critical group.
(The edge-probability must be bounded away from $1$ to avoid
complete graphs and complete multipartite graphs, but I don't
really know what the `right' bound should be.)
There is some experimental evidence for this, but it is not
extensive.  The following weak form is probably more accessible.
\begin{CONJ}
If $c>1$ and $c\log(n)/n<p(n)<1-o(\log(n)/n)$ then, as
$n\goesto\infty$, the expected value of $\nu(\G(n,p))$ 
remains bounded.
\end{CONJ}
By considering Smith normal forms and using Proposition 2.3,
it is easy to see that for a connected undirected graph $G$,
if $\kappa(G)$ is square-free then $K(G)$ is cyclic (or trivial).
Since the density of square-free natural numbers is asymptotically
$6/\pi^{2}$, this motivates the third conjecture.
\begin{CONJ}
If $c>1$ and $c\log(n)/n<p(n)<1-o(\log(n)/n)$ then, as
$n\goesto\infty$, the probability that $\kappa(\G(n,p))$ is
square-free is $(1-o(1))6/\pi^{2}$.
\end{CONJ}

\section{Equitable partitions of graphs.}

See Chapter 5 of Godsil \cite{Go} for further development and 
application of the theory of equitable partitions of undirected
graphs.

Consider a graph $G$, and let
$\pi=(\pi_{1},\ldots,\pi_{p})$ be an ordered partition of $V$ into
pairwise disjoint nonempty blocks.  The partition $\pi$ is
\emph{equitable} for $G$ provided that there exist nonnegative
integers $F_{ij}$ and $R_{ij}$ for all $1\leq i,j\leq p$ such that every
vertex in $\pi_{i}$ is the initial vertex of exactly $F_{ij}$ directed
edges of $G$ which have their terminal vertices in $\pi_{j}$, and
every vertex in $\pi_{j}$ is the terminal vertex of exactly $R_{ij}$
directed edges of $G$ which have their initial vertices in $\pi_{i}$.
(The letters $F$ and $R$ are mnemonic for `forward' and `reverse',
respectively.)  These integers define $p\times p$ matrices $F$ and $R$,
and we regard $F$ as the adjacency matrix of a graph
$G/\pi$ on the vertex-set $\{1,\ldots,p\}$, called the \emph{quotient
of $G$ by $\pi$}.  It will be convenient to use the notations $A$,
$\Delta$ and $Q$ for $A(G)$, $\Delta(G)$, and $Q(G)$, and to use
$F$, $D$ and $\widehat{Q}$ for $A(G/\pi)$, $\Delta(G/\pi)$, and $Q(G/\pi)$.

For $\pi$ an equitable partition of $G$, let $P$ be the matrix
indexed by $V\times\{1,\ldots,p\}$, with entries
$$P_{vi}:=\left\{\begin{array}{ll}
1 &\ \mathrm{if}\ v\in\pi_{i},\\
0 &\ \mathrm{if}\ v\not\in\pi_{i}.\end{array}\right.$$
Then $B:=P^{\dagger}P$ is the invertible $p\times p$ diagonal matrix with
entries $B_{ii}=\#\pi_{i}$ for each $1\leq i\leq p$.  For each
$1\leq i,j\leq p$, by counting in two ways the directed edges of $G$ with initial
vertex in $\pi_{i}$ and terminal vertex in $\pi_{j}$ we see that
$B_{ii}F_{ij}=R_{ij}B_{jj},$
yielding the matrix equations $BF=RB$ and $B^{-1}R=FB^{-1}$.
For any vertex $v\in\pi_{i}$ we have
$$\Delta_{vv}=\sum_{w\in V}A_{vw}=\sum_{j=1}^{p}F_{ij}=D_{ii},$$
or, in matrix form, $\Delta=PDB^{-1}P^{\dagger}$.  Therefore, $\Delta P=PD$.
Also, for $v\in\pi_{i}$, and any $1\leq j\leq p$,
$$(AP)_{vj}=\sum_{w\in\pi_{j}}A_{vw}=F_{ij}=(PF)_{vj},$$
so that $AP=PF$.  It follows that $QP=P\widehat{Q}$.
Finally, consider any $1\leq i\leq p$ and $v\in\pi_{j}$.  Then
$$(B^{-1}P^{\dagger}A)_{iv}=B_{ii}^{-1}\sum_{w\in\pi_{i}}A_{wv}=B_{ii}^{-1}
R_{ij}=F_{ij}B_{jj}^{-1}=(FB^{-1}P^{\dagger})_{iv},$$
so that $B^{-1}P^{\dagger}A=FB^{-1}P^{\dagger}$.  Also, since
$$B^{-1}P^{\dagger}\Delta=B^{-1}P^{\dagger}PDB^{-1}P^{\dagger}=DB^{-1}P^{\dagger},$$
we conclude that $B^{-1}P^{\dagger}Q=\widehat{Q}B^{-1}P^{\dagger}$.

\section{Critical groups of graph quotients.}

We continue with the notation of the previous section.  The matrix
$P$ defines a group homomorphism $P:\ZZ^{p}\goesto\ZZ^{V}$, and
$P\widehat{Q}\ZZ^{p}=QP\ZZ^{V}\subseteq Q\ZZ^{V}$.
Therefore, $P$ induces a homomorphism
$$\rho:\K^{*}(G/\pi)\longrightarrow\K^{*}(G)$$
between the dual critical groups, which is
well-defined by $\rho(\x+\widehat{Q}\ZZ^{p}):=P\x+Q\ZZ^{V}.$

Theorem 6.1 was inspired by Theorem 10.2 of Biggs \cite{Bi3}.
\begin{THM}  Let $G$ be a strongly connected graph, and
let $\pi$ be an equitable partition of $G$.  Then the natural
homomorphism $\rho:\K^{*}(G/\pi)\goesto\K^{*}(G)$ is injective.
\end{THM}
\begin{proof}
To show that $\rho:\K^{*}(G/\pi)\goesto\K^{*}(G)$ is injective, we must
show that if $\x\in\ZZ^{p}$ is such that $P\x\in Q\ZZ^{V}$, then
$\x\in \widehat{Q}\ZZ^{p}$.  Accordingly, assume that $\x\in\ZZ^{p}$
and $\v\in\ZZ^{V}$ satisfy $P\x=Q\v$.  Let $\y:=B^{-1}P^{\dagger}\v$.  Then
$$\x=B^{-1}P^{\dagger}P\x=B^{-1}P^{\dagger}Q\v=\widehat{Q}B^{-1}P^{\dagger}\v=
\widehat{Q}\y.$$
The entries of $\y$ are rational numbers, but we need to find 
$\u\in\ZZ^{p}$ such that $\x=\widehat{Q}\u$.  To do this, notice that
$$Q\v=P\x=P\widehat{Q}\y=QP\y,$$
so that $Q(\v-P\y)=\zero$.  Since $G$ is strongly connected,
the kernel of $Q$ is $\RR\one$, and this implies that $\v-P\y=c\one$ for some $c\in\RR$.  Now
$\v=P\y+c\one=P(\y+c\one)$, and since every entry of $\v$ is an
integer, every entry of $\u:=\y+c\one$ is an integer.  Finally,
since $\widehat{Q}\one=\zero$ it follows that $\x=\widehat{Q}\y=
\widehat{Q}\u$, which shows that $\x\in\widehat{Q}\ZZ^{p}$ and
completes the proof.
\end{proof}
Of course, under the hypotheses of Theorem 6.1, the natural homomorphism
$$\rho:K^{*}(G/\pi)\longrightarrow K^{*}(G)$$
is also injective.

\begin{EG}\emph{
The hypothesis that $G$ is strongly connected can not be dropped
from Theorem 6.1, as the following example shows.  Let
$$Q=\left[\begin{array}{rrr} 2 & -1 & -1\\ 0 & 0 & 0\\ 0 & 0 & 0
\end{array}\right]\hspace{1cm}\mathrm{and}\hspace{1cm}
\widehat{Q}=\left[\begin{array}{rr} 2 & -2 \\ 0 & 0
\end{array}\right].$$
Then $Q$ is the Laplacian matrix of a weakly connected graph $G$ with
vertex-set $\{1,2,3\}$ and $\widehat{Q}$ is the Laplacian of the 
quotient of $G$ by the equitable partition $\pi=\{\{1\},\{2,3\}\}$ of
$G$.  However, by computing the Smith normal forms of $Q$ and
$\widehat{Q}$ one sees that
$$\K^{*}(G)\simeq\ZZ\oplus\ZZ\hspace{1cm}\mathrm{and}\hspace{1cm}
\K^{*}(G/\pi)\simeq(\ZZ/2\ZZ)\oplus\ZZ.$$
Hence, the homomorphism $\rho:\K^{*}(G/\pi)\goesto\K^{*}(G)$ is not
injective in this case.
}\end{EG}

\begin{EG}\emph{
With the hypotheses of Theorem 6.1, $K^{*}(G/\pi)$ is regarded as
a subgroup of $K^{*}(G)$ \emph{via} the natural inclusion $\rho$.
However, $K^{*}(G/\pi)$ might not be a direct summand of
$K^{*}(G)$, as the following example shows.  The nine-cycle
$C_{9}$ has an equitable partition $\pi$ for which the
quotient graph $C_{9}/\pi$ is the three-cycle $C_{3}$,
but $K^{*}(C_{3})\simeq\ZZ/3\ZZ$ is not a direct summand
of $K^{*}(C_{9})\simeq\ZZ/9\ZZ$.
}\end{EG}

\section{Matroid invariance of $K(G)$.}

Let $G$ be a strongly connected graph, and let $\h$ be the
vector of vertex activities defined in Proposition 3.2.  If vertex
$v\in V$ is such that $h(v)=1$ then we say that $v$ is
a \emph{simple} vertex of $G$.

\begin{LMA}  Let $G$ be a strongly connected graph, and let $v,w\in
V(G)$.  If $w$ is simple then the Laplacian matrix $Q(G)$ is
equivalent to the matrix $Q'$ obtained from $Q$ by
replacing every entry in either column $v$ or row $w$ by zero.
\end{LMA}
\begin{proof}
Let $\h$ be the vector of vertex activities of $G$.
Use elementary column operations to add column $u$ to column $v$ for
all $v\neq u\in V$.  Then use elementary row operations to add $h(u)$ 
times row $u$ to row $w$ for all $w\neq u\in V$.  The resulting
matrix $Q'$ is equivalent to $Q$ and has the required form, since
$Q^{\dagger}\h=Q\one=\zero$.
\end{proof}

\begin{PROP}  Let $G$ and $H$ be vertex-disjoint strongly connected
graphs, and let $v\in V(G)$ and $w\in V(H)$.  Denote by
$(G\cup H)/vw$ the graph obtained from $G\cup H$ by identifying
$v$ and $w$.  If $w$ is a simple vertex of $H$ then
$$\K(G\cup H)\simeq \K((G\cup H)/vw)\oplus\ZZ$$
and
$$K(G\cup H)\simeq K((G\cup H)/vw).$$
\end{PROP}
\begin{proof}
Let $M$ be the submatrix of $Q(G)$ obtained by deleting the
row and column indexed by $v$; so we have
$$Q(G)=\left[\begin{array}{ll}
M & \overline{\v} \\ \v & \Delta(G)_{vv}\end{array}\right]$$
for some row vector $\v$ and column vector $\overline{\v}$.
Since $Q(G)\one=\zero$, we have
$$Q(G)\approx Q'(G):=\left[\begin{array}{ll}
M & \zero \\ \v & 0\end{array}\right].$$
Similarly, let $N$ be the submatrix of $Q(H)$ obtained
by deleting the row and column indexed by $w$; so we have
$$Q(H)=\left[\begin{array}{ll}
\Delta(H)_{ww} & \w \\ \overline{\w} & N\end{array}\right]$$
for some row vector $\w$ and column vector $\overline{\w}$.
Since $w$ is a simple vertex of $H$, Lemma 7.1 implies that
$Q(H)\approx[0]\oplus N$.

Now, since $Q(G\cup H)=Q(G)\oplus Q(H)$ we see that
$Q(G\cup H)\approx Q'(G)\oplus[0]\oplus N$.  Also,
since $w$ is a simple vertex of $H$ we see that
$$Q((G\cup H)/vw)=\left[\begin{array}{lll}
M & \overline{\v} & O \\
\v & d & \w \\
O & \overline{\w} & N \end{array}\right]\approx
\left[\begin{array}{lll}
M & \zero & O \\
\v & 0 & \zero \\
O & \zero & N \end{array}\right]
=Q'(G)\oplus N$$
(here, $d:=\Delta(G)_{vv}+\Delta(H)_{ww}$).
Since $Q(G\cup H)\approx
Q((G\cup H)/vw)\oplus[0]$ we have $\snf(G\cup H)=
\snf((G\cup H)/vw) \oplus[0]$, from which the result follows.
\end{proof}

Let $G$ and $H$ be vertex-disjoint weakly connected
directed graphs, let $v_{1}\neq v_{2}$ be distinct vertices of $G$, and
let $w_{1}\neq w_{2}$ be distinct vertices of $H$.  The graph
$G\bullet H:=(G\cup H)/\{v_{1}w_{1},v_{2}w_{2}\}$ is obtained from
$G\cup H$ by identifying $v_{1}$ and $w_{1}$, and identifying $v_{2}$
and $w_{2}$.  The graph $G\circ H:=(G\cup H)/\{v_{1}w_{2},v_{2}w_{1}\}$
is obtained from $G\cup H$ by identifying $v_{1}$ and $w_{2}$, and
identifying $v_{2}$ and $w_{1}$.  We say that $G\bullet H$ and 
$G\circ H$
are related by \emph{twisting a two-vertex cut}.

\begin{PROP}  Let $G$ and $H$ be vertex-disjoint strongly connected
graphs, and use the notation of the above paragraph.  If both $w_{1}$
and $w_{2}$ are simple vertices of $H$ then
$$\K(G\bullet H)\simeq\K(G\circ H).$$
\end{PROP}
\begin{proof}
For $i=1,2$, let $\v_{i}$ be the row vector with entries
$Q(G)_{v_{i}z}$ for each $z\in V(G)\drop\{v_{1},v_{2}\}$,
and let $\overline{\v}_{i}$ be the column vector with entries
$Q(G)_{zv_{i}}$ for each $z\in V(G)\drop\{v_{1},v_{2}\}$.
For $i=1,2$, let $\w_{i}$ be the row vector with entries
$Q(H)_{w_{i}z}$ for each $z\in V(H)\drop\{w_{1},w_{2}\}$,
and let $\overline{\w}_{i}$ be the column vector with entries
$Q(H)_{zw_{i}}$ for each $z\in V(H)\drop\{w_{1},w_{2}\}$.
Then the matrices $Q(G\bullet H)$ and $Q(G\circ H)$ have the block
forms as shown: 

$$\left[\begin{array}{llll}
M & \overline{\v}_{1} & \overline{\v}_{2} & O \\
\v_{1} & a_{11} & a_{12} & \w_{1}\\
\v_{2} & a_{21} & a_{22} & \w_{2}\\
O & \overline{\w}_{1} & \overline{\w}_{2} & N
\end{array}\right]
\hspace{1cm} \mathrm{and} \hspace{1cm}
\left[\begin{array}{llll}
M & \overline{\v}_{1} & \overline{\v}_{2} & O \\
\v_{1} & b_{11} & b_{12} & \w_{2}\\
\v_{2} & b_{21} & b_{22} & \w_{1}\\
O & \overline{\w}_{2} & \overline{\w}_{1} & N
\end{array}\right]$$
$$Q(G\bullet H)\hspace{5cm} Q(G\circ H)$$
Here, $M$ is the submatrix of $Q(G)$
induced by rows and columns in $V(G)\drop\{v_{1},v_{2}\}$,
$N$ is the submatrix of $Q(H)$ induced by rows and columns in $V(H)
\drop\{w_{1},w_{2}\}$, and 
$$a_{ij}:=Q(G)_{v_{i}v_{j}}+Q(H)_{w_{i}w_{j}}\ \ \mathrm{and}\ \
b_{ij}:=Q(G)_{v_{i}v_{j}}+Q(H)_{w_{3-i},w_{3-j}}$$
for $1\leq i\leq 2$.

Since both $w_{1}$ and $w_{2}$ are simple vertices of $H$, 
an argument analogous to the proof of Lemma 7.1 shows that
$Q(G\bullet H)$ and $Q(G\circ H)$ are equivalent
to
$$\left[\begin{array}{llll}
M & \overline{\v}_{1} & \zero & O \\
\v_{1} & a_{11} & 0 & \w_{1}\\
\v_{1}+\v_{2} & c_{21} & 0 & \zero \\
O & \overline{\w}_{1} & \zero & N
\end{array}\right]
\hspace{1cm} \mathrm{and} \hspace{1cm}
\left[\begin{array}{llll}
M & \overline{\v}_{1} & \zero & O \\
\v_{1} & b_{11} & 0 & \w_{2}\\
\v_{1}+\v_{2} & c_{21} & 0 & \zero\\
O & \overline{\w}_{2} & \zero & N
\end{array}\right]$$
$$Q'(G\bullet H)\hspace{5cm} Q'(G\circ H)$$
respectively, in which $c_{21}=Q(G)_{v_{1}v_{1}}+Q(G)_{v_{2}v_{1}}$.

Let $\h$ denote the vector of activities of $H$, so that 
$\h^{\dagger}Q(H)=\zero$.  Then the sum over all $z\in V(H)\drop\{w_{1},
w_{2}\}$ of $h(z)$ times the $z$-th row of $N$ is equal to
$-\w_{1}-\w_{2}$, since both $w_{1}$ and $w_{2}$ are simple in $H$.
Since the columns of $Q(H)$ sum to $\zero$,
the columns of $N$ sum to $-\overline{\w}_{1}-\overline{\w}_{2}$.
Now, for $z\in V(H)\drop\{w_{1},w_{2}\}$, add $h(z)$ times row $z$
of $Q'(G\circ H)$ to the row indexed by $v_{1}w_{2}$.  Then,
for $z\in V(H)\drop\{w_{1},w_{2}\}$, add column $z$ of the
resulting matrix to the column indexed by $v_{1}w_{2}$.
The result is the matrix
$$Q''(G\circ H)=
\left[\begin{array}{rrrr}
M & \overline{\v}_{1} & \zero & O \\
\v_{1} & c_{11} & 0 & -\w_{1}\\
\v_{1}+\v_{2} & c_{21} & 0 & \zero\\
O & -\overline{\w}_{1} & \zero & N
\end{array}\right]$$
in which
\showon
c_{11}
&=& b_{11}+\sum_{z\in V(H)\drop\{w_{1},w_{2}\}}h(z)\overline{w}_{2}(z)
-\sum_{z\in V(H)\drop\{w_{1},w_{2}\}}w_{1}(z)\\
&=& b_{11}-Q(H)_{w_{2}w_{2}}-Q(H)_{w_{2}w_{1}}+
Q(H)_{w_{1}w_{1}}+Q(H)_{w_{2}w_{1}}\\
&=&  b_{11}-\Delta(H)_{w_{2}w_{2}}+\Delta(H)_{w_{1}w_{1}}\\
&=& \Delta(G)_{v_{1}v_{1}}+\Delta(H)_{w_{2}w_{2}}
-\Delta(H)_{w_{2}w_{2}}+\Delta(H)_{w_{1}w_{1}}\\
&=& a_{11}
\showoff
Finally, by multiplying the last $n(H)-2$ rows and columns of
$Q''(G\circ H)$ by $-1$, we obtain the matrix $Q'(G\bullet H)$.
This shows that $Q(G\bullet H)$ and $Q(G\circ H)$ are equivalent, so
that $\snf(Q(G\bullet H))=\snf(Q(G\circ H))$ and $\K(G\bullet H)\simeq
\K(G\circ H)$.
\end{proof}

\begin{EG}\emph{
The hypothesis in Proposition 7.3 that $w_{1}$ and $w_{2}$ are simple
vertices of $H$ can not be removed, as the following example shows.
Let $G$ and $H$ both have Laplacian matrix
$$Q(G)=Q(H)=\left[\begin{array}{rr} 1 & -1 \\ -2 & 
2\end{array}\right],$$
and note that the vertex activities are $2$ and $1$ in each graph.
Then 
$$Q(G\bullet H)=\left[\begin{array}{rr} 2 & -2 \\ -4 & 
4\end{array}\right]\hspace{1cm}\mathrm{and}\hspace{1cm}
Q(G\circ H)=\left[\begin{array}{rr} 3 & -3 \\ -3 & 
3\end{array}\right],$$
and by calculating Smith normal forms we see that
$$\K(G\bullet H)\simeq(\ZZ/2\ZZ)\oplus\ZZ
\hspace{1cm}\mathrm{and}\hspace{1cm}
\K(G\circ H)\simeq(\ZZ/3\ZZ)\oplus\ZZ.$$
}\end{EG}

A strongly connected directed graph $G$ is \emph{balanced} if,
for each vertex $v\in V(G)$, the indegree of $v$ equals the outdegree
of $v$.  Undirected graphs are thus a special case of balanced graphs.
Equivalently, $G$ is balanced if and only if $Q^{\dagger}\one=\zero$,
\emph{i.e.} every vertex of $G$ is simple.

\begin{CORO}  Let $G$ and $H$ be vertex-disjoint balanced graphs,
and let $v\in V(G)$ and $w\in V(H)$.  Then
$$\K(G\cup H)\simeq \K((G\cup H)/vw)\oplus\ZZ$$
and
$$K(G\cup H)\simeq K((G\cup H)/vw).$$
\end{CORO}

\begin{CORO}  Let $G$ and $H$ be vertex-disjoint balanced
graphs, let $v_{1}\neq v_{2}$ in $V(G)$, and let $w_{1}\neq w_{2}$
in $V(H)$.  With the notation introduced above,
$$\K(G\bullet H)\simeq\K(G\circ H).$$
\end{CORO}

Whitney \cite{Wh} characterized those pairs of undirected graphs
$G, H$ which have isomorphic graphic matroids as being
exactly those pairs for which $G$ may be transformed into $H$
by some sequence of splittings or mergings of one-vertex cuts
and twistings of two-vertex cuts.  Corollary 7.7 follows immediately.

\begin{CORO}  Let $G$ and $H$ be undirected graphs. If the graphic
matroids of $G$ and $H$ are isomorphic, then $K(G)\simeq K(H)$.
\end{CORO}

\begin{EG}\emph{
The converse of Corollary 7.7 does not hold.  In fact, for undirected
graphs, the Tutte polynomial $T(G;x,y)$ is not computable from the
critical group $\K(G)$, as the following example shows.  Let
$v$ be a vertex of $C_{3}$ and let $w$ be a vertex of $C_{4}$, where
the cycles are vertex-disjoint.  Then
\showon
\K((C_{3}\cup C_{4})/vw) &\simeq& 
(\ZZ/3\ZZ)\oplus(\ZZ/4\ZZ)\oplus\ZZ\\
&\simeq &(\ZZ/12\ZZ)\oplus\ZZ\ \simeq\ \K(C_{12}).
\showoff
However,
$$T((C_{3}\cup C_{4})/vw;x,y)=(y+x+x^{2})(y+x+x^{2}+x^{3})$$
and
$$T(C_{12};x,y)=y+x+x^{2}+\cdots+x^{11}.$$
}\end{EG}

Finally, I conjecture that the critical group is not
computable from the Tutte polynomial, either.
\begin{CONJ} There exist connected undirected graphs
$G$ and $H$ such that $T(G;x,y)=T(H;x,y)$ and 
$\K(G)\not\simeq \K(H)$.
\end{CONJ}

\section{Inequivalence of $\T(G)$ and $K(G)$.}

For this section we consider only connected undirected graphs.
By Proposition 2.3, in this case the reduced critical
group $K(G)$ is a finite abelian group of order $\kappa(G)$, the
cardinality of the set $\T(G)$ of spanning trees of $G$.  It is thus
reasonable to consider the problem of constructing a bijection from
$\T(G)$ to $K(G)$ in a `natural' way.   We will show that, in general,
this is not possible.  To allow more flexibility in the construction,
we consider the complex vector spaces $\CC\T(G)$
and $\CC K(G)$, and ask for an isomorphism of vector spaces
$\psi_{G}:\CC\T(G)\goesto\CC K(G)$ which is constructed naturally
from $G$.  The naturality condition means that $\psi_{G}$ should
depend only on the isomorphism class of $G$, but we must first state
this more precisely.

For any graph isomorphism $f:G\goesto H$, there is an induced
bijection $f_{\T}:\T(G)\goesto\T(H)$ defined by sending each spanning
tree of $G$ to its image under $f$.  This extends linearly to an
isomorphism from $\CC\T(G)$ to $\CC\T(H)$ which we also denote by
$f_{\T}$.
Let $\Phi$ be the matrix indexed by $V(H)\times V(G)$, with
$$\Phi_{vw}:=\left\{\begin{array}{ll}
1 &\ \mathrm{if}\ f(w)=v,\\
0 &\ \mathrm{if}\ f(w)\neq v.\end{array}\right.$$
Then $\Phi:\ZZ^{V(G)}\goesto\ZZ^{V(H)}$ is an isomorphism, and
$\Phi Q^{\dagger}(G)=Q^{\dagger}(H)\Phi$.  It follows that $\Phi 
Q^{\dagger}(G)\ZZ^{V(G)}=
Q^{\dagger}(H)\ZZ^{V(H)}$, and so $f$ induces a group isomorphism
$$f_{\K}:\K(G)\stackrel{\sim}{\longrightarrow}\K(H)$$
well-defined by $f_{\K}(\x+Q^{\dagger}(G)\ZZ^{V(G)}):=
\Phi\x+Q^{\dagger}(H)\ZZ^{V(H)}.$
For reduced critical groups, we also have a group isomorphism
$$f_{K}:K(G)\stackrel{\sim}{\longrightarrow}K(H)$$
induced by $f$, and we extend this linearly to obtain an
isomorphism from $\CC K(G)$ to $\CC K(H)$, also denoted by $f_{K}$.

The naturality condition on $\psi_{G}$ is that, for any multigraph
isomorphism $f:G\goesto H$, the diagram
$$\begin{array}{lcr}
\CC\T(G)  & \stackrel{\psi_{G}}{\longrightarrow} & \CC K(G) \\
{\scriptstyle f_{\T}}\ \downarrow & &
\downarrow \ {\scriptstyle f_{K}}\\
\CC\T(H) & \stackrel{\psi_{H}}{\longrightarrow} & \CC K(H)\end{array}$$
is commutative; that is to say, $\psi_{H}\circ f_{\T}=f_{K}\circ\psi_{G}$.
This means that $\psi_{H}$ is `the same as' $\psi_{G}$, after 
relabelling the vertices according to $f$.

\begin{THM}  There exist connected undirected graphs  $G$ for
which there is no natural isomorphism $\psi_{G}:\CC\T(G)\goesto
\CC K(G)$.
\end{THM}
\begin{proof}
In particular, the naturality condition must hold when $G=H$ and
$f$ is in the automorphism group $\aut(G)$ of $G$.  In this situation,
the assignment $f\mapsto f_{\T}$ gives a representation of $\aut(G)$
acting on $\CC\T(G)$ and the assignment $f\mapsto f_{K}$ gives a
representation of $\aut(G)$ acting on $\CC K(G)$.  (All the
representation theory we need is in Chapter 1 of Ledermann \cite{L}.)
Commutativity of
the diagram means that $\psi_{G}$ is an $\aut(G)$-equivariant 
isomorphism, so these two representations are linearly equivalent.
A representation of a finite group is determined up to linear
equivalence by its group character, so $\psi_{G}$ exists if and only
if the characters $\chi_{\T}$ and $\chi_{K}$ of these two 
representations are equal.  Since these are permutation 
representations, their characters are given by counting fixed points;
that is, for each $f\in\aut(G)$,
$$\chi_{\T}(f)=\#\{T\in\T(G):\ f_{\T}(T)=T\}$$
and
$$\chi_{K}(f)=\#\{\x\in K(G):\ f_{K}(\x)=\x\},$$
respectively.

Thus, to show that a natural construction of $\psi_{G}$ is impossible,
it suffices to find a connected undirected graph $G$ and automorphism
$f\in\aut(G)$ such that $\chi_{\T}(f)\neq\chi_{K}(f)$.  This is easy:
let $G$ be a circulant graph with at least three vertices,
and let $f\in\aut(G)$ be a cyclic permutation of all of $V(G)$.
Every spanning tree of $G$ has both leaves and non-leaf vertices, and
so it can not be left fixed by $f$.  Thus, $\chi_{\T}(f)=0$.  On the other
hand, the $\zero$ element of $K(G)$ is such that $f_{K}(\zero)=\zero$,
since $f_{K}$ is a group automorphism.  Thus, $\chi_{K}(f)\geq 1$,
and since $\chi_{\T}\neq\chi_{K}$ it follows that $\psi_{G}$ does
not exist.
\end{proof}

I have not thought much about the following problem, but the
question seems interesting.
\begin{PROB}\emph{
With notation as in the proof of Theorem 8.1,
does there exist a vertex-transitive connected undirected
graph $G$, with at least three vertices, such that the characters
$\chi_{\T}$ and $\chi_{K}$ of $\aut(G)$ are equal?
}\end{PROB}
Note, in particular, that for such a graph every automorphism
fixes at least one spanning tree;  together with vertex-transitivity,
this seems to be quite restrictive.

In contrast with Theorem 8.1, a negative result, if one
fixes a total order
$\prec$ on the edge-set of $G$ then the automorphism group of the
pair $(G,\prec)$ is trivial (as long as $G$ has at least three
vertices), so the naturality condition for these structures becomes
vacuous and the possibility arises of constructing a bijection from
$\T(G)$ to $K(G)$ relative to $\prec$.  Jonathan Dumas
(personal communication, August 2000) has recently found
such a construction which, moreover, behaves well regarding
the external activities of spanning trees with respect to $\prec$.

\section{The dollar game for strongly connected graphs.}

For undirected graphs, a combinatorial understanding of the critical
group has been developed by Biggs \cite{Bi1,Bi2,Bi3}.  We generalize
this to the case of strongly connected graphs;  the theory is 
essentially the same as in the undirected case, with one interesting
extra complication.

Let $G=(V,E)$ be a strongly connected graph with no loops, and let
$\bank$ denote a designated vertex of $G$, which we call the \emph{bank}.
A \emph{configuration} is an integer vector $\c\in\ZZ^{V}$ such that
$\one^{\dagger}\c=0$, so that $c(\bank)=-\sum_{v\neq\bank}c(v)$.
We say that a configuration $\c$ is \emph{nonnegative} when
$c(v)\geq 0$ for all $\bank\neq v\in V$.  A vertex
$v\neq\bank$ is \emph{legal} for $\c$ when $c(v)\geq
\Delta_{vv}$, the outdegree of $v$.  The configuration $\c$ is 
\emph{stable} when $c(v)<\Delta_{vv}$ for all 
$v\neq\bank$; in this case, and in this case only, the bank vertex
$\bank$ is
\emph{legal} for $\c$.  For a configuration $\c$ and vertex $v$, the effect 
of \emph{firing $v$ from} $\c$ is the configuration 
$\c|v:=\c-Q^{\dagger}\hat{v}$, in which $\hat{v}\in\ZZ^{V}$ denotes the 
characteristic vector of the vertex $v\in V$.  More explicitly,
for each $w\in V$, $(c|v)(w)$ is defined by
$$(c|v)(w)=\left\{\begin{array}{ll}
c(v)-\Delta_{vv} & \mathrm{if}\ w=v,\\
c(w)+A_{vw} & \mathrm{if}\ w\neq v.\end{array}\right.$$
If one considers a nonnegative configuration $\c$ as representing
$c(v)$ dollars at each vertex $v\neq\bank$, then the effect of firing
a legal vertex $v$ is to send one dollar along each edge with initial
vertex $v$.  More generally, if $v_{1}\cdots v_{k}$ is a sequence of
vertices in which $v\in V$ occurs $m(v)$ times, then the effect of
firing this sequence of vertices from the configuration $\c$ is
$\c|v_{1}\cdots v_{k}=\c-Q^{\dagger}\m$, in which 
$\m=\hat{v}_{1}+\cdots+\hat{v}_{k}$ is the vector
of multiplicities.  The sequence $v_{1}\cdots v_{k}$ is
\emph{legal for} $\c$ provided that $v_{i}$ is legal for 
$\c|v_{1}\cdots v_{i-1}$ for each $1\leq i\leq k$.
For a configuration $\c$, let $\S(\c)$ denote the set of all 
configurations $\b$ such that $\b=\c|v_{1}\cdots v_{k}$
for some sequence $v_{1}\cdots v_{k}$ of vertices which is legal for $\c$
and does not contain the bank vertex $\bank$.

\begin{PROP}  Let $G=(V,E)$ be a loopless strongly connected graph with
bank vertex $\bank$. For every configuration $\c$ on
$(G,\bank)$, the set $\S(\c)$ is finite.  If $\c$ is nonnegative then
every configuration in $\S(\c)$ is also nonnegative.
\end{PROP}
\begin{proof}
If $\c$ is nonnegative and $v$ is legal for $\c$, then $\c|v$ is
nonnegative.  From this it follows that if $\c$ is nonnegative
then every configuration in $\S(\c)$ is nonnegative.  More
generally, for any configuration $\c$, define the
configuration $\c^{-}$ by
$$c^{-}(v):=\left\{\begin{array}{ll}
c(v) & \mathrm{if}\ v\neq\bank\ \mathrm{and}\ c(v)<0,\\
0 & \mathrm{if}\ v\neq\bank\ \mathrm{and}\ c(v)\geq 0,
\end{array}\right.$$
and $c^{-}(\bank):=-\sum_{v\neq\bank}c^{-}(v)$.  If $v_{1}\cdots
v_{k}$ is a legal sequence for $\c$ which does not contain $\bank$,
then it is a legal sequence for $\c-\c^{-}$.  Since $\c-\c^{-}$ is
nonnegative, it follows that $\c|v_{1}\cdots v_{k}-\c^{-}$ is
nonnegative.  That is, $\b-\c^{-}$ is nonnegative for all $\b\in
S(\c)$.

For each $v\in V$, let $d(v)$ denote the length of a shortest directed
path from $v$ to $\bank$ in $G$, and let $r:=\max\{d(v):\ v\in V\}$.
For each configuration $\c$, define the `label' of $\c$ to be
$\ell(\c):=(\ell_{1},\ldots,\ell_{r})$, in which $\ell_{i}:=
\sum\{c(v):\ v\in V\ \mathrm{and}\ d(v)=i\}$.  Define a total order
$\prec$ on $\ZZ^{r}$ as follows:\ $(p_{1},\ldots,p_{r})\prec
(q_{1},\ldots,q_{r})$ if and only if either $p_{1}+\cdots+p_{r}<
q_{1}+\cdots+q_{r}$, or $p_{1}+\cdots+p_{r}=q_{1}+\cdots+q_{r}$ and
$p_{1}=q_{1}$, $p_{2}=q_{2}$,\ldots $p_{i-1}=q_{i-1}$, $p_{i}>q_{i}$
for some $1\leq i\leq r$.  Notice that $(\ZZ^{r},\prec)$ has the
same order type as $(\ZZ,<)$.  Also notice that if $v\neq\bank$
and $v$ is legal for $\c$, then
$$\ell(\c^{-})\preceq\ell((\c|v)^{-})\preceq 
\ell(\c|v)\prec\ell(\c).$$
It follows that, for any configuration $\c$, the set $\{\ell(\b):\
\b\in\S(\c)\}$ is finite.  But for any $\q\in\ZZ^{r}$, the set
of configurations $\{\b:\ \ell(\b)=\q\}$ is also finite.  These
two observations suffice to show that $\S(\c)$ is finite.
\end{proof}

\begin{LMA} Let $G=(V,E)$ be a loopless strongly connected graph with
bank vertex $\bank$.  Let $\c$ be a configuration on $(G,\bank)$,
let $v,w\in V$, and let $\c|w:=\c-Q^{\dagger}\hat{w}$.  If $v\neq w$
then $(c|w)(v)\geq c(v)$.  In particular, if $v\not\in\{w,\bank\}$ and
$v$ is legal for $\c$, then $v$ is legal for $\c|w$.
\end{LMA}
\begin{proof}
This follows immediately from the facts that the off-diagonal elements
of $Q$ are nonpositive, and that a vertex $v$ is legal for a 
configuration $\c$ if and only if $c(v)\geq\Delta_{vv}$.
\end{proof} 

\begin{LMA}  Let $G=(V,E)$ be a loopless strongly connected graph with
bank vertex $\bank$.  Let $\c$ be a configuration on $(G,\bank)$,
let $\m\in\NN^{V}$, and let $v_{1}\cdots v_{k}$ be a sequence of
vertices, with multiplicity vector $\n:=\hat{v}_{1}+\cdots
+\hat{v}_{k}$.  Produce the sequence $w_{1}\cdots w_{\ell}$ by 
deleting the first $\min\{m(z), n(z)\}$ occurrences of vertex $z$ from 
the sequence $v_{1}\cdots v_{k}$, for each $z\in V$.  If $v_{1}\cdots
v_{k}$ is legal for $\c$ and $n(\bank)=0$, then $w_{1}\cdots w_{\ell}$
is legal for $\c':=\c-Q^{\dagger}\m$.
\end{LMA}
\begin{proof}
We proceed by induction on $k$, the length of $v_{1}\cdots v_{k}$.

For the basis of induction, $k=1$, assume that $v\neq\bank$ is
legal for $\c$.  If $m(v)\geq 1$ then the empty sequence is 
legal for $\c'$, as required.  Otherwise, write
$\m=\hat{u}_{1}+\cdots+\hat{u}_{r}$ for some sequence of vertices
$u_{1}\cdots u_{r}$.  Since $m(v)=0$, $v$ does not occur in the
sequence $u_{1}\cdots u_{r}$.  The previous lemma and induction
on $r$ now show that $v$ is legal for $\c'$, as required.

For the induction step, assume the result for sequences of
length $k-1$, and consider $v_{1}\cdots v_{k}$.  First,
assume that $m(v_{1})\geq 1$, and consider $\b:=\c|v_{1}$ and
$\m':=m-\hat{v}_{1}$.  Then $\c'=\b-Q^{\dagger}\m'$ and
$v_{2}\cdots v_{k}$ is legal for $\b$.  Applying the induction
hypothesis to $\b$, $\m'$, and $v_{2}\cdots v_{k}$, we see that
$w_{1}\cdots w_{\ell}$ is legal for $\c'$.  For the remaining
case, assume that $m(v_{1})=0$, so that $w_{1}=v_{1}$.  As in the
basis of induction, since $v_{1}$ is legal for $\c$, $v_{1}$ is
legal for $\c'$.  Now apply the induction hypothesis to
$\b:=\c|v_{1}$, $\m$, and $v_{2}\cdots v_{k}$ to conclude that
$w_{2}\cdots w_{\ell}$ is legal for $\c'|v_{1}$.  Hence,
$w_{1}\cdots w_{\ell}$ is legal for $\c'$, completing the
induction step and the proof.
\end{proof}

\begin{PROP}  Let $G=(V,E)$ be a loopless strongly connected graph with
bank vertex $\bank$.  For every  configuration $\c$
on $(G,\bank)$, the set $\S(\c)$ contains a unique stable
configuration.
\end{PROP}
\begin{proof}
Let $\D$ be the graph with vertex-set $\S(\c)$ and directed edges
$\b\goesto\b|v$ when $v\neq\bank$ is legal for $\b\in\S(\c)$.
Then $\D$ is a nonempty graph, and, by the proof of Proposition 9.1,
since $\ell(\b|v)\prec\ell(\b)$ for all $\b\in\S(\c)$
$\D$ contains no directed cycles.  Therefore, $\D$ has at least
one sink vertex, which is a stable configuration on $(G,\bank)$.

Now suppose that $\a$ and $\b$ are two stable configurations
in $\S(\c)$.  Let $v_{1}\cdots v_{k}$ and $u_{1}\cdots u_{r}$
be sequences of vertices which are legal for $\c$, do not
contain $\bank$, and are such that $\c|v_{1}\cdots v_{k}=\a$
and $\c|u_{1}\cdots u_{r}=\b$.  Let $\n:=\hat{v}_{1}+\cdots+
\hat{v}_{k}$ and let $\m:=\hat{u}_{1}+\cdots+\hat{u}_{r}$.
The hypothesis of Lemma 9.3 is satisfied, so produce
the subsequence $w_{1}\cdots w_{\ell}$ of $v_{1}\cdots v_{k}$
as in that lemma.  Now, since $w_{1}\cdots w_{\ell}$ is a 
legal sequence for $\c-Q^{\dagger}\m=\c|u_{1}\cdots u_{r}=\b$
which does not contain $\bank$, and since $\b$ is stable,
it follows that $w_{1}\cdots w_{\ell}$ is the empty sequence.
From the construction of $w_{1}\cdots w_{\ell}$, it follows
that $m(z)\geq n(z)$ for each $z\in V$.  By symmetry, we may
repeat this argument with $v_{1}\cdots v_{k}$ and $u_{1}\cdots
u_{r}$ interchanged, and deduce that $n(z)\geq m(z)$ for each
$z\in V$.  Finally, since $\m=\n$ we conclude that
$\a=\c-Q^{\dagger}\n=\c-Q^{\dagger}\m=\b$, finishing
the proof.
\end{proof}

We define the \emph{stabilization} of a
configuration $\c$ to be the unique stable configuration in
$\S(\c)$.  If $\c$ is a stable configuration, then we define
the \emph{successor} $\sigma(\c)$ of $\c$
to be the stabilization of $\c|\bank$.  Thus, $\sigma$ is
an endofunction on the set of stable
configurations on $(G,\bank)$.  We say that a 
stable configuration $\c$ is \emph{critical} when
$\sigma^{k}(\c)=\c$ for some positive integer $k$.
(Here, $\sigma^{k}$ denotes the $m$-th functional iterate
of $\sigma$.)

\begin{LMA}
Let $G=(V,E)$ be a loopless strongly connected graph with
bank vertex $\bank$.\\
\noindent\textup{(a)}\  Every critical configuration
on $(G,\bank)$ is nonnegative.\\
\noindent\textup{(b)}\  For every stable configuration
$\c$ on $(G,\bank)$, there is a nonnegative integer
$m$ such that $\sigma^{m}(\c)$ is nonnegative.\\
\noindent\textup{(b)}\  For every stable configuration
$\c$ on $(G,\bank)$, there is a nonnegative integer
$m$ such that $\sigma^{m}(\c)$ is critical.
\end{LMA}
\begin{proof}
For part (a), let $\c$ be a critical configuration, and
let $v_{1}\cdots v_{k}$ be a nonempty legal sequence of
vertices for $\c$, such that $\c|v_{1}\cdots v_{k}=\c$.
Then $\n:=\hat{v}_{1}+\cdots+\hat{v}_{k}$ is a nonzero
vector of nonnegative integers such that $\c=\c-Q^{\dagger}\n$,
so that $Q^{\dagger}\n=\zero$.
By Proposition 3.2, it follows that $\n=\lambda\h$ is a
positive integer multiple of the vector $\h$ of vertex
activities. Every $v\in V$ is fired in the sequence
$v_{1}\cdots v_{k}$ exactly $\lambda h(v)$ times, hence
at least once, and so $c(v)\geq 0$ for all $v\neq\bank$.

For part (b), we use the strategy of the proof of
Proposition 9.1, but with a different definition of the label
of a configuration.  For a configuration
$\c$, define $\c^{-}$ as in the proof of Proposition 9.1.  If 
$\c^{-}\neq\zero$, then for each $v\neq\bank$, let $d(v)$
denote the length of a shortest directed path which begins
at $v$ and ends at a vertex $w\neq\bank$ such that $c(w)<0$,
and let $r$ be the maximum value of $d(v)$ for all
$v\neq\bank$ such that $c(v)\geq 0$.  For
$1\leq j\leq r$, let $\ell_{j}:=\sum\{c(v):\ d(v)=j\}$,
and define the `label' of $\c$ to be 
$\ell(\c):=(\c^{-},\ell_{1},\ldots,\ell_{r})$.
Define a partial order on the set of all labels as follows:
$(\a,k_{1},\ldots,k_{s})\prec(\b,\ell_{1},\ldots,\ell_{r})$
if either $\b\neq\a$ and $\b-\a$ is nonnegative (except at $\bank$), or
$\b=\a\neq\zero$ (in which case $s=r>0$) and $k_{1}=\ell_{1}$,
$k_{2}=\ell_{2}$,\ldots, $k_{j-1}=\ell_{j-1}$, $k_{j}<\ell_{j}$ for some
$1\leq j\leq r$.

Notice that for any configuration $\a$ with $-\a$ nonnegative,
there are only finitely many stable
configurations $\c$ such that $\c^{-}=\a$ (in
fact, the number of them is $\prod\{\Delta_{vv}:\ v\neq\
\bank\ \mathrm{and}\ a(v)=0\}$).  It follows that the
set of all stable configurations, partially ordered by
the order $\prec$ on their labels, has no infinite
ascending chains.  Now, if $v$ is legal for
$\c$ and $\c^{-}\neq\zero$, then $\ell(\c)\prec\ell(\c|v)$,
even when $v=\bank$.  Hence, if $\c$ is stable but
not nonnegative, then $\ell(\c)\prec\ell(\sigma(\c))$.
Since there are no infinite ascending chains, there
is a nonnegative integer $m$ such that $\sigma^{m}(\c)$
is maximal, and this must be a stable, nonnegative
configuration.

For part (c), let $\c$ be any stable configuration.
By part (b), we may assume that $\c$ is nonnegative.
There are exactly $\prod\{\Delta_{vv}: v\neq\bank\}$
nonnegative stable configurations on $(G,\bank)$,
and the successor function
acts as an endofunction on this set.  Since this
set is finite, for any $\c$ in it there exists a nonnegative
integer $m$ and positive integer $k$ such that
$\sigma^{m+k}(\c)=\sigma^{m}(\c)$.  That is,
$\sigma^{m}(\c)$ is critical.
\end{proof}

Let $\C(G,\bank)$ denote the set of all
critical configurations of $(G,\bank)$.  For critical
configurations $\a$ and $\b$ on $(G,\bank)$, we say that
$\a$ and $\b$ are \emph{coeval} if $\sigma^{m}(\a)=\b$
for some nonnegative integer $m$.

\begin{LMA}
Let $G=(V,E)$ be a loopless strongly connected graph with
bank vertex $\bank$.   Coevalence is an equivalence
relation on $\C(G,\bank)$, and each Coevalence class
has cardinality divisible by $h(\bank)$, in which
$\h$ is the vector of vertex activities.  
\end{LMA}
\begin{proof}
That coevalence is an equivalence relation is easy to see.
Let $\c$ be a critical configuration on $(G,\bank)$, let
$v_{1}\cdots v_{k}$ be a nonempty sequence of vertices
which is legal for $\c$, such that $\c|v_{1}\cdots v_{k}=\c$,
and as short as possible subject to these conditions,
and let $\n:=\hat{v}_{1}+\cdots+\hat{v}_{k}$.  As in the
proof of Lemma 9.5(a), we see that $\n=\lambda\h$ for
some positive integer $\lambda$.  Since the bank
vertex occurs $\lambda h(\bank)$ times in the sequence
$v_{1}\cdots v_{k}$ it follows that $\c$ is coeval
with exactly $\lambda h(\bank)$ critical configurations.
\end{proof}

\begin{PROP}  Let $G=(V,E)$ be a loopless strongly connected graph,
and let $\bank$ be a simple vertex of $G$.  Then 
every coevalence class of $\C(G,\bank)$ is a singleton.
\end{PROP}
\begin{proof}
Let $\h$ be the vector of vertex activities of $G$.
Let $\c$ be any critical configuration of $(G,\bank)$, and 
let $v_{1}\cdots v_{k}$ be a nonempty sequence of vertices
which is legal for $\c$ and such that $\c|v_{1}\cdots v_{k}=\c$.
Let $\n:=\hat{v}_{1}+\cdots+\hat{v}_{k}$ be the multiplicity
vector of $v_{1}\cdots v_{k}$; as in Lemma 9.5(a)
we have $\n=\lambda\h$ for some positive integer $\lambda$.

Since $\c$ is stable, $v_{1}=\bank$.  If this is the
only occurrence of $\bank$ in the sequence $v_{1}\cdots v_{k}$,
then $\sigma(\c)=\c$ and it follows that the coevalence class
of $\c$ is the singleton $\{\c\}$.  Otherwise,
let $1=p_{1}<p_{2}<\cdots<p_{r}\leq k$ be all the indices
from $1$ to $k$ such that $v_{p_{i}}=\bank$, and let $p_{r+1}=k+1$.
Since $\n=\lambda\h$ and $h(\bank)=1$, we see that $\lambda= r$.

Now, $v_{2}\cdots v_{p_{2}-1}$ is a legal
sequence for $\c|\bank$ which does not contain any occurrence
of $\bank$.  Applying Lemma 9.3 with $\m=\h-\hat{\bank}$, let
$w_{1}\cdots w_{\ell}$ be the sequence so produced.  This
sequence is legal for $(\c|\bank)-Q^{\dagger}\m=\c-Q^{\dagger}\h
=\c$, and since $\c$ is stable, it follows that $w_{1}\cdots w_{\ell}$
is the empty sequence.  Therefore, for each $\bank\neq u\in V$,
the multiplicity of $u$ in $v_{2}\cdots v_{p_{2}-1}$ is at most
$h(u)$. 

Since $\b:=\c|\bank v_{2}\cdots v_{p_{2}-1}$ is critical, we
may repeat the argument of the above paragraph, using $\b|\bank$
and $v_{p_{2}+1}\cdots v_{p_{3}-1}$.  As above, we conclude that
for each $\bank\neq u\in V$, the multiplicity of $u$ in
$v_{p_{2}+1}\cdots v_{p_{3}-1}$ is at most $h(u)$. 
For each $1\leq i\leq r$, let $\n_{i}:=\hat{v}_{p_{i}+1}+\cdots
+\hat{v}_{p_{i+1}-1}$.  Applying the above argument for each
sequence $v_{p_{i}+1}\cdots v_{p_{i+1}-1}$, we see that
$n_{i}(u)\leq h(u)$ for each $1\leq i\leq r$ and $\bank\neq u\in V$.
But $\n_{1}+\cdots+\n_{r}+r\hat{\bank}=\n=r\h$, from which it follows
that $\n_{1}=\cdots=\n_{r}=\h-\hat{\bank}$.  Therefore,
$\c|\bank v_{2}\cdots v_{p_{2}-1}=\c-Q^{\dagger}\h=\c$, so that
$\sigma(\c)=\c$,  and it follows that the coevalence class
of $\c$ is the singleton $\{\c\}$, completing the proof.
\end{proof}

\begin{EG}\emph{
As the following example shows, if $\bank$ is not a 
simple vertex of $G$, then it may happen that some
coevalence class of $\C(G,\bank)$ has cardinality
strictly greater than $h(\bank)$.  Consider the
graph $G$ with Laplacian matrix and vector of vertex
activities
$$Q=\left[\begin{array}{rrrr}
2 & -1 & -1 & 0\\
0 & 2 & -2 & 0\\
-1 & 0 & 2 & -1\\
-1 & -1 & 0 & 2\end{array}\right]
\hspace{1cm}\mathrm{and}\hspace{1cm}
\h=\left[\begin{array}{r}
3\\ 5\\ 8\\ 4\end{array}\right].$$
Denoting the vertices by $1$, $2$, $3$, $4$ corresponding to
the matrix indices, take $\bank=1$ for the bank vertex.  Denoting
a configuration $\c$ for this $(G,\bank)$ by the triple $c(2)c(3)c(4)$,
we see that the nonnegative stable configurations are mapped by
the successor function as
$$\begin{array}{ccccccccccccc}
000 & & & & 101 & & & & & & & &\\
\downarrow & & & & \downarrow & & & & & & & &\\
110 & \goesto & 100 & \goesto & 011 & \goesto & 001 & \goesto & 111 & 
\goesto & 101 & \goesto & 110.
\end{array}$$
Thus, $\C(G,\bank)$ has a single coevalence class of cardinality six.
}\end{EG}

Finally, we connect the critical group of $G$ with the results
of this section;  this should be compared with Theorems 3.8
and 8.1 of Biggs \cite{Bi3}.
Let $\c$ be any configuration on $(G,\bank)$, and let $\overline{\c}$
be the stabilization of $\c$.  From Lemma 9.5(c), there is a
nonnegative integer $m_{0}$ such that $\sigma^{m}(\overline{\c})$
is critical for all $m\geq m_{0}$.  These critical
configurations are all coeval with one another, and this
coevalence class is determined uniquely by $\c$;  we denote
it by $[\c]$.

\begin{THM}   Let $G=(V,E)$ be a loopless strongly connected graph with
bank vertex $\bank$, and let $U\subseteq\ZZ^{V}$ consist of those
vectors $\u\in\ZZ^{V}$ such that $\one^{\dagger}\u=0$.  Then $K(G)$
is the subgroup $U/Q^{\dagger}\ZZ^{V}$ of $\K(G)$, and the
elements of $K(G)$ correspond bijectively with coevalence classes 
of critical configurations on $(G,\bank)$ \emph{via} the
correspondence $\u+Q^{\dagger}\ZZ^{V}\mapsto [\u]$.
\end{THM}
\begin{proof}
For $\x\in\ZZ^{V}$ the element $\z+Q^{\dagger}\ZZ^{V}$ of $\K(G)$ is in $K(G)$
if and only if $\lambda\z\in Q^{\dagger}\ZZ^{V}$ for some positive 
integer $\lambda$.  Since every column $\q$ of $Q^{\dagger}$ satisfies
$\one^{\dagger}\q=0$, if $\z+Q^{\dagger}\ZZ^{V}$ is in $K(G)$ then
$\one^{\dagger}\z=0$.  Conversely, if $\one^{\dagger}\z=0$ then
$\z=Q^{\dagger}\w$ for some rational $V$-indexed vector $\w$, since
the rank of $Q^{\dagger}$ is $n(G)-1$.  Hence, $\lambda\z\in 
Q^{\dagger}\ZZ^{V}$ for some positive integer $\lambda$, so that
$\z+Q^{\dagger}\ZZ^{V}$ is in $K(G)$.  This proves that 
$K(G)=U/Q^{\dagger}\ZZ^{V}$.

For the second claim, consider two configurations $\b$ and $\c$ on
$(G,\bank)$, so that $\b,\c\in U$.  If $[\b]=[\c]$ then there are
nonnegative integers $p,q$ such that $\sigma^{p}(\overline{\b})=\sigma^{q}
(\overline{\c})$.  Therefore, there are $V$-indexed vectors
$\m,\n$ of nonnegative integers such that $\b-Q^{\dagger}\m=
\c-Q^{\dagger}\n$, so that $\b-\c$ is in $Q^{\dagger}\ZZ^{V}$.
This proves that the map $\u+Q^{\dagger}\ZZ^{V}\mapsto[\u]$ is
injective.  Since this map is clearly surjective, the theorem is 
proved.
\end{proof}

\begin{PROB}\emph{
Given a strongly connected graph $G$ and bank vertex $\bank\in V$, say that
$\bank$ is \emph{small} when every coevalence class of $\C(G,\bank)$ has
cardinality $h(\bank)$, and that $\bank$ is \emph{fair} when all 
coevalence classes of $\C(G,v)$ have equal cardinality.  Are there
polynomial-time algorithms for determining whether a given vertex
is small, or is fair?
}\end{PROB}

\end{document}